\documentclass[11pt,oneside]{amsart}

\usepackage{microtype}

\usepackage{amsmath,ifthen, amsfonts, amssymb,
srcltx, amsopn, color, enumerate} 
\usepackage[linktocpage=true]{hyperref}
\usepackage[cmtip,arrow]{xy}
\usepackage{pb-diagram, pb-xy}
\usepackage{overpic}

 \def\notorth{\:\ensuremath{\reflectbox{\rotatebox[origin=c]{90}{$\nvdash$}}}}

\newcommand{\showcomments}{yes}

\newsavebox{\commentbox}

\newcounter{ax}
\newtheorem{thm}{Theorem}[section]
\newtheorem{lem}[thm]{Lemma}

\theoremstyle{definition}
\newtheorem{defn}[thm]{Definition}
\newtheorem{rem}[thm]{Remark}

\newtheorem{claim*}{Claim}

\newtheorem{cons}[thm]{Construction}

\DeclareMathOperator{\dimension}{dim}

\DeclareMathOperator{\diam}{\textup{\textsf{diam}}}

\DeclareMathOperator{\hull}{hull}

\newcommand{\neb}{\mathcal N}

\newcommand{\field}[1]{\mathbb{#1}}

\newcommand{\naturals}{\ensuremath{\field{N}}}
\newcommand{\reals}{\ensuremath{\field{R}}}

\newcommand{\HHS}[2]{\ensuremath{(\cuco {#1},\frak {#2})}}

\newcommand{\boundary}{{\ensuremath \partial}}

\makeatletter

\newcommand{\Rmnum}[1]{\mathbf{{\expandafter\@slowromancap\romannumeral #1@}}}

\makeatother

\newcommand{\tup}[1]{\vec{#1}}

\let\oldmarginpar\marginpar
\renewcommand\marginpar[1]{\-\oldmarginpar[\raggedleft\footnotesize #1]{\raggedright\footnotesize #1}}

\newcommand{\tsh}[1]{\left\{\kern-.7ex\left\{#1\right\}\kern-.7ex\right\}}
\newcommand{\Tsh}[2]{\tsh{#2}_{#1}}
\newcommand{\ignore}[2]{\Tsh{#2}{#1}}

\newcommand{\co}{\colon}

\newcounter{enumitemp}

\newcommand{\dist}{\textup{\textsf{d}}}

\newcommand{\cuco}[1]{{\mathcal #1}}

\newcommand{\fontact}{{\mathcal C}}

\newcommand{\gate}{\mathfrak g}

\newcommand{\coneoff}[1]{\widehat{#1}}

\newcommand{\nest}{\sqsubseteq}
 \usepackage{mathabx}
\newcommand{\propnest}{\sqsubsetneq}
\newcommand{\orth}{\bot}
\newcommand{\transverse}{\pitchfork}

\newcommand{\median}{\mathfrak m}

\begin{document}

\title{What is a hierarchically hyperbolic space?}

\author[A. Sisto]{Alessandro Sisto}
\address{Department of Mathematics, ETH Zurich, 8092 Zurich, Switzerland}
\email{sisto@math.ethz.ch}

\maketitle

\begin{abstract}
 The first part of this survey is a heuristic, non-technical discussion of what an HHS is, and the aim is to provide a good mental picture both to those actively doing research on HHSs and to those who only seek a basic understanding out of pure curiosity. It can be read independently of the second part, which is a detailed technical discussion of the axioms and the main tools to deal with HHSs.
\end{abstract}

\setcounter{tocdepth}{2}
\setcounter{secnumdepth}{2}
\tableofcontents

\section*{Introduction}

Hierarchically hyperbolic spaces (HHSs) were introduced in \cite{HHS1} as a common framework to study mapping class groups and cubical groups. The definition is inspired by the extremely successful Masur-Minsky machinery to study mapping class groups \cite{MM1,MM2,Behrstock_thesis}. Since \cite{HHS1}, the list of examples has expanded significantly \cite{HHS2,HHS_asdim,HagenSusse}, and the HHS framework has been used to prove several new results, including new results for mapping class groups and cubical groups. For example, the previously known bound for the asymptotic dimension of mapping class groups has been dramatically improved in \cite{HHS_asdim}, while the main result from \cite{HHS_quasiflats} is that top-dimensional quasi-flats in HHSs stay within bounded distance from a finite union of ``standard orthants'', a fact that was known neither for mapping class groups nor for cubical groups without imposing additional constraints (see \cite{Huang:quasiflats}).

The aim of this survey article, however, is not to present the state of the art of the field, which is very much evolving. In this direction, we only give a brief description of all the relevant papers below. The main aim of this survey is, instead, to discuss the geometry of HHSs, only assuming that the reader is familiar with (Gromov-)hyperbolic spaces. The definition of HHS is admittedly hard to digest if one is not presented with the geometric intuition behind it, and the aim of this survey is to remedy this shortcoming. The first part is aimed at the casual reader and gives a general idea of what an HHS looks like. We will discuss the various notions in the main motivating examples too; the reader can use whichever example they are familiar with to gain better understanding.

The second part of the survey is mainly aimed at those who want to do research on HHSs, as well as to those who seek deeper understanding. We will discuss every axiom in detail, and then we will proceed to discuss the main tools one can use to study HHSs. Anyone who becomes familiar with the material that will be presented will have a rather deep understanding of HHSs. And will be ready to tackle one of the many open questions asked in the papers we describe below...

\subsection*{State of the art}

\begin{itemize}
 \item In \cite{HHS1}, J. Behrstock, M. Hagen and I axiomatised the Masur-Minsky machinery, extended it to right-angled Artin groups (and many other groups acting on CAT(0) cube complexes including fundamental groups of special cube complexes), and initiated the study of the geometry of hierarchically hyperbolic groups by studying quasi-flats via ``coarse differentiation''.
 \item In \cite{HHS2} we simplified the list of axioms, which allowed us to extend the list of examples of hierarchically hyperbolic groups (and to significantly simplify the Masur-Minsky approach). The paper contains a combination theorem for trees of HHSs, as well as other results to construct new HHSs out of old ones. 
 \item Speaking of new examples, in \cite{HagenSusse} M.Hagen and T. Susse prove that all proper cocompact CAT(0) cube complexes are HHSs.
 \item \cite{HHS_asdim} deals with asymptotic dimension. We show finiteness of the asymptotic dimension of hierarchically hyperbolic groups, giving explicit estimates in certain cases. In the process we drastically improve previously known bounds on the asymptotic dimension of mapping class groups. We also show that many natural (small cancellation) quotients of hierarchically hyperbolic groups are hierarchically hyperbolic.
 \item In \cite{HHS_boundary}, M. Durham, M. Hagen and I introduce a compactification of hierarchically hyperbolic groups, related to Thurston's compactification of Teichm\"uller space in the case of mapping class groups and Teichm\"uller space. This compactification turns out to be very well-behaved as, for example, ``quasi-convex subgroups'' in a suitable sense have a well-defined and easily recognisable limit set (while the situation for Teichm\"uller space is more complicated). Constructing this compactification allowed us, for example, to study dynamical properties of individual elements and to prove a rank rigidity result.
 \item Further study of the HHS boundary is carried out in \cite{Mousley1,Mousley2}, where S. Mousley shows non-existence of boundary maps in certain cases and other exotic phenomena.
 \item  In \cite{HHS_quasiflats}, J. Behrstock, M. Hagen and I study the geometry of quasi-flats, that is to say images of quasi-isometric embeddings of $\mathbb  R^n$ in hierarchically hyperbolic spaces. More specifically, we show that top-dimensional quasi-flats lie within finite distance of a union of ``standard orthants''. This simultaneously solves open questions and conjectures for most of the motivating examples of hierarchically hyperbolic groups, for example a conjecture of B. Farb for mapping class groups, and one by J. Brock for the Weil-Petersson metric, and it is new in the context of CAT(0) cube complexes too.
 \item In \cite{Spriano}, D. Spriano studies non-trivial HHS structures on hyperbolic spaces, and uses them to show that certain natural amalgamated products of hierarchically hyperbolic groups are hierarchically hyperbolic.
 \item In \cite{HHS_stable}, C. Abbott, J. Behrstock and M. Durham prove that hierarchically hyperbolic groups admit a ``best'' acylindrical action on a hyperbolic space, and provide a complete classification of stable subgroups of hierarchically hyperbolic groups.
 \item In \cite{Haettel:morphism}, T. Haettel studies homomorphisms of higher rank lattices to hierarchically hyperbolic groups, finding severe restrictions.
 \item In \cite{largest_random}, using part of an HHS structure, S. Taylor and I studied various notions of projections, including subsurface projection for mapping class groups, along a random walk, and used this to prove a conjecture of I. Rivin on random mapping tori.
\end{itemize}

\subsection*{Acknowledgements} This article has been written for the proceedings for the ``Beyond hyperbolicity'' conference held in June 2016 in Cambridge, UK. The author would like to thank Mark Hagen, Richard Webb, and Henry Wilton for organising the conference and the wonderful time he had in Cambridge.

The author would also like to thank Jason Behrstock, Mark Hagen, and Davide Spriano for useful comments on previous drafts of this survey.

\part{Heuristic discussion}\label{part:heuristic}

\section{Standard product regions}
In this section we discuss the first heuristic picture of an HHS, which is the one provided by \emph{standard product regions}.

If an HHS $\cuco X$ is not hyperbolic, then the obstruction to its hyperbolicity is encoded by the collection of its standard product regions. These are quasi-isometrically embedded subspaces that split as direct products, and the crucial fact is that each standard product region, as well as each of its factors, is an HHS itself, and in fact an HHS of lower ``complexity''. It is not very important at this point, but the complexity is roughly speaking the length of a longest chain of standard product regions $P_1\subsetneq P_2\subsetneq \dots\subsetneq P_n$ contained in the HHS; what is important right now is that factors of standard product regions are ``simpler'' HHSs, and the ``simplest'' HHSs are hyperbolic spaces. This is what allows for induction arguments, where the base case is that of hyperbolic spaces.

Standard product regions encode entirely the non-hyperbolicity of the HHS $\cuco X$ in the following sense. Given a, say, length metric space $(Z,d)$ and a collection of subspaces $\mathcal P$, one can define the cone-off of $Z$ with respect to the collection of subspaces (in several different ways that coincide up to quasi-isometry, for example) by setting $d'(x,y):=1$ for all $x,y$ contained in the same $P\in\mathcal P$ and $d'(x,y)=d(x,y)$ otherwise, and declaring the cone-off distance between two points $x,y$ to be $\inf_{x=x_0,\dots,x_n=y}\sum d'(x_i,x_{i+1})$. This has the effect of collapsing all $P\in\mathcal P$ to bounded sets, and the reason why this is a sensible thing to do is that one might want to consider the geometry of $Z$ ``up to'' the geometry of the $P\in\mathcal P$. When $Z$ is a graph, as is most often the case for us, coning-off amounts to adding edges connecting pairs of vertices contained in the same $P\in\mathcal P$.

Back to HHSs, when coning-off all standard product regions of an HHS one obtains a hyperbolic space, that we denote $\fontact S$\footnote{This notation is taken from the mapping class group context, even though it's admittedly not the best notation in other examples.}. In other words, an HHS is weakly hyperbolic relative to the standard product regions. Roughly speaking, when moving around $\cuco X$, one is either moving in the hyperbolic space $\fontact S$ or in one of the standard product regions. The philosophy behind many induction arguments for HHSs is that when studying a certain ``phenomenon'', either it leaves a visible trace in $\fontact S$, or it is ``confined'' in a standard product region, and can hence be studied there. For example, if the HHS is in fact a group, one can consider the subgroup generated by an element $g$, and it turns out that either the orbit maps of $g$ in $\fontact S$ are quasi-isometric embeddings, or $g$ virtually fixes a standard product region \cite{HHS_boundary}.

So far we discussed the ``top-down'' point of view on standard product regions, but there is also a ``bottom-up'' approach. In fact, one can regard HHSs as built up inductively starting from hyperbolic spaces, in the following way:
\begin{itemize}
 \item hyperbolic spaces are HHSs,
 \item direct products of HHSs are HHSs,
 \item ``hyperbolic-like'' arrangements of HHSs are HHSs.
\end{itemize}

The third bullet refers to $\fontact S$ being hyperbolic, and the fact that $\fontact S$ can also be thought of as encoding the intersection pattern of standard product regions. Incidentally, I believe that there should be a characterisation of HHSs that looks like the list above, i.e. that by suitably formalising the third bullet one can obtain a characterisation of HHSs. This has not been done yet, though. There is, however, a combination theorem for trees of HHSs in this spirit \cite{HHS2}.

One final thing to mention is that standard product regions have well-behaved coarse intersections, meaning that the coarse intersection of two standard product regions is well-defined and coarsely coincides with some standard product region. In other words, $\cuco X$ is obtained gluing together standard product regions along sub-HHSs, so a better version of the third bullet above would be ``hyperbolic-like arrangements of HHS glued along sub-HHSs are HHSs''.

\subsection{In the examples}

We now discuss standard product regions in motivating examples of HHSs.

\subsubsection{RAAGs}
Consider a simplicial graph $\Gamma$. Whenever one has a (full) subgraph $\Lambda$ of $\Gamma$ which is the join of two (full, non-empty) subgraphs $\Gamma_1,\Gamma_2$, then the RAAG $A_\Gamma$ contains an undistorted copy of the RAAG $A_\Lambda\approx A_{\Lambda_1}\times A_{\Lambda_2}$. Such subgroups and their cosets are the standard product regions of $A_\Gamma$. In this case, $\fontact S$ is a Cayley graph of $A_\Gamma$ with respect to an infinite generating set (unless $\Gamma$ consists of a single vertex), namely $V\Gamma\cup\{A_\Lambda<A_\Gamma:\Lambda=join(\Lambda_1,\Lambda_2)\}$. A given HHS can be given different HHS structures (which turns out to allow for more flexibility when performing various constructions, rather than being a drawback), and one instance of this is that one can regard as standard product regions all $A_\Lambda$ where $\Lambda$ is any proper subgraph of $\Gamma$, one of the factors being trivial. In this case $\fontact S$ is the Cayley graph of $A_\Gamma$ with respect to the generating set $V\Gamma\cup \{A_\Lambda<A_\Gamma:\Lambda{\rm\ proper\ subgraph\ of\ }\Gamma\}$, which is perhaps more natural.
%To re-iterate the philosophy explained above, if one can study a given problem about $A_\Gamma$ in some proper sub-RAAG, it is best to do so, otherwise $\fontact S$ probably provides insight. 

For both HHS structures described above, $\fontact S$ is not only hyperbolic, but in fact quasi-isometric to a tree.

\subsubsection{Mapping class groups}
Given a surface $S$, there are some ``obvious'' subgroups of $MCG(S)$ that are direct products. In fact, consider two disjoint (essential) subsurfaces $Y,Z$ of $S$. Any two self-homeomorphisms of $S$ supported respectively on $Y$ and $Z$ commute. This yields (up to ignoring issues related to the difference between boundary components and punctures that I do not want to get into) a subgroup of $MCG(S)$ isomorphic to $MCG(Y)\times MCG(Z)$. Such subgroups are in fact undistorted. One can similarly consider finitely many disjoint subsurfaces instead, and this yields the standard product regions in $MCG(S)$. More precisely, one should fix representatives of the (finitely many) topological types of collections of disjoint subsurfaces, and consider the cosets of the subgroups as above. In terms of the marking graph, product regions are given by all markings containing a given sub-marking.

In this case it shouldn't be too hard to convince oneself that $\fontact S$ as defined above is quasi-isometric to the curve complex, see \cite[]{MM1}. To re-iterate the philosophy explained above, if some behaviour within $MCG(S)$ is not confined to a proper subsurface $Y$, then the geometry of $\fontact S$ probably comes into play when studying it, and otherwise it is most convenient to study the problem on the simpler subsurface $Y$.

\subsubsection{CAT(0) cube complexes}

Hyperplanes are crucial for studying CAT(0) cube complexes, and the carrier of a hyperplane (meaning the union of all cubes that the hyperplane goes through) is naturally a product of the hyperplane and an interval. It is then natural, when trying to define an HHS structure on a CAT(0) cube complex, to regard carriers of hyperplanes as standard product regions, even though one of the factors is bounded.

This is not enough, though. As mentioned above, coarse intersections of standard product regions should be standard product regions. But it is easy to describe the coarse intersection of two carriers of hyperplanes, or more generally the coarse intersection of two convex subcomplexes. Given a convex subcomplex $Y$ of the CAT(0) cube complex $\cuco X$, one can consider the \emph{gate map} $\gate_Y:\cuco X\to Y$, which is the closest-point projection in either the CAT(0) or the $\ell^1$--metric (they coincide). More combinatorially, for $x\in\cuco X^{(0)}$, $\gate_Y(x)$ is defined by the property that the hyperplanes separating $x$ from $\gate_Y(x)$ are exactly those separating $x$ from $Y$. The coarse intersection of convex subcomplexes $Y,Z$ is just $\gate_Y(Z)$, which is itself a convex subcomplex.

Back to constructing an HHS structure on a CAT(0) cube complex, we now know that we need to include as standard product regions all gates of carriers in other carriers. But then we are not done yet, because for the same reason that we need to take gates of carriers we also need to take gates of gates, and so on. Also, we need this process to stabilise eventually (which is not always the case, unfortunately), because an HHS needs to have finite complexity to allow for induction arguments. All these considerations lead to the definition of \emph{factor system}. Rather than carrier of hyperplanes, we will consider combinatorial hyperplanes, which are the two copies of a hyperplane that bound its carrier, but this does not make a substantial difference for the purposes of this discussion.

\begin{defn}
 A factor system $\mathcal F$ for the cube complex $\cuco X$ is a collection of convex subcomplexes so that:
 \begin{enumerate}
  \item all combinatorial hyperplanes are in $\mathcal F$,
  \item there exists $\xi\geq 0$ so that if $F,F'\in\mathcal F$ and $\gate_{F}(F')$ has diameter at least $\xi$, then $\gate_{F}(F')\in\mathcal F$.
  \item $\mathcal F$ is uniformly locally finite.
 \end{enumerate}
\end{defn}

Any factor system $\mathcal F$ on a cube complex gives an HHS structure, where $\fontact S$ is obtained coning off all members of $\mathcal F$. It turns out that $\fontact S$ is quasi-isometric to a tree. It is proven in \cite{HagenSusse} that all cube complexes admitting a proper cocompact action by isometries have a factor system, and are therefore HHSs.

\section{Projections to hyperbolic spaces}
In this section we discuss a point of view on HHSs that is more similar to the actual definition, and it is in terms of ``coordinates'' in certain hyperbolic spaces.

We already saw that any HHS $\cuco X$ comes equipped with a hyperbolic space, $\fontact S$, obtained collapsing the standard product regions. In particular, there is a coarsely Lipschitz map $\pi_S:\cuco X\to \fontact S$.

The (coarse geometry of the) hyperbolic space $\fontact S$ is not enough to recover the whole geometry of $\cuco X$, since it does not contain information about the standard product regions themselves. Hence, we want something to keep track of the geometry of the standard product regions. Since factors of standard product regions are HHSs themselves, they also come with a hyperbolic space obtained collapsing the standard product sub-regions. Considering all standard product regions, we obtain a collection of hyperbolic space $\{\fontact Y\}_{Y\in\mathfrak S}$, which, together, control the geometry of $\cuco X$, as we are about to discuss. The index set $\mathfrak S$ is the set of factors of standard product regions, where the whole of $\cuco X$ should be considered as a product region with (a trivial factor and the other factor being) $S\in\mathfrak S$, so as to include $\fontact S$ among the hyperbolic spaces we consider.

Another piece of data we need is a collection of coarsely Lipschitz maps $\pi_Y:\cuco X\to\fontact Y$ for all $Y\in\mathfrak S$, which allow us to talk about the geometry of $\cuco X$ ``from the point of view of $\fontact Y$''. These projection maps come from natural coarse retractions of $\cuco X$ onto the standard product regions, composed with the collapsing maps, but for now it only matters that the $\pi_Y$ exist. 

\subsection{Distance formula and hierarchy paths}
The first way in which the $\fontact Y$ control the geometry of $\cuco X$ is that whenever $x,y\in\cuco X$ are far away, then their projections to some $\fontact Y$ are far away, so that any coarse geometry feature of $\cuco X$ leaves a trace in at least one of the $\fontact Y$. In fact, there is much better control on distances in $\cuco X$ in terms of distances in the various $\fontact Y$, and this is given by the \emph{distance formula}. This is perhaps the most important piece of machinery in the HHS world, and certainly the most iconic. To state it, we need a little bit of notation. We write $A\approx_K B$ if $A/K-K\leq B\leq KB+K$, and declare $\ignore{A}{L}=A$ if $A\geq L$, and $\ignore{A}{L}=0$ otherwise. The distance formula says that for all sufficiently large $L$ there exists $K$ so that
$$\dist_{\cuco X}(x,y)\approx_K \sum_{Y\in\mathfrak S} \ignore{\dist_{\fontact Y}(\pi_Y(x),\pi_Y(y))}{L}$$
for all $x,y\in\cuco X$. In words, the distance in $\cuco X$ between two points is, up to multiplicative and additive constants, the sum of the distances between their far-away projections in the various $\fontact Y$. Very imprecisely, this is saying that $\cuco X$ quasi-isometrically embeds in $\prod_{Y\in\mathfrak S} \fontact Y$ endowed with some sort of $\ell^1$ metric. To save notation one usually writes $\dist_Y(x,y)$ instead of $\dist_{\fontact Y}(\pi_Y(x),\pi_Y(y))$.

Another important fact related to the distance formula is the existence of \emph{hierarchy paths}, that is to say quasigeodesics in $\cuco X$ that shadow geodesics in each $\fontact Y$. Namely, there exists $D$ so that for any $x,y\in\cuco X$ there exists a $(D,D)$--quasigeodesic $\gamma$ joining them so that $\pi_Y\circ\gamma$ is an unparametrised $(D,D)$--quasigeodesic in $\fontact Y$. Since $\fontact Y$ is hyperbolic, being an unparametrised quasigeodesic means that (the image of) $\pi_Y\circ\gamma$ is Hausdorff-close to a geodesic, and it ``traverses'' the geodesic coarsely monotonically. In most cases it is much better to deal with hierarchy paths than with geodesics.

\subsection{Consistency}
The distance formula alone is not enough for almost anything, but the point is that it comes with a toolbox that one uses to control the various projection terms by constraining certain projections in terms of certain other projections. In Part \ref{part:technical}, we will analyse these tools in detail. For now we will instead describe what happens to various projections when moving along a hierarchy path, which gives the right picture about how projections are constrained.

The tl; dr version of this subsection is: for certain pairs $Y,Z\in\mathfrak S$, along a hierarchy path one can only change the projections to $Y,Z$ in a specified order. This is sufficient to understand most of the next subsection.

\subsubsection{Nesting} Let $x,y\in\cuco X$ and suppose that $\dist_Y(x,y)$ is large. Notice that $Y$ (which is a factor of a standard product region) gives a bounded set in $\fontact S$ (which was obtained from $\cuco X$ by collapsing standard product regions), which we denote $\rho^Y_S$. We know that when moving along any hierarchy path $\gamma$ from $x$ to $y$, the projection to $Y$ needs to change. This is how this happens: $\gamma$ has an initial subpath where the projection to $Y$ coarsely does not change, while the projection to $\fontact S$ approaches $\rho^Y_S$. All the progress that needs to be made by $\gamma$ in $\fontact Y$ is made by a middle subpath whose projection to $\fontact S$ remains close to $\rho^Y_S$. Then, there is a final subpath that does not make any progress in $\fontact Y$ and takes us from $\rho^Y_S$ to $\pi_S(y)$ in $\fontact S$. In short, you can only make progress in $\fontact Y$ if you are close to $\rho^Y_S$ in $\fontact S$.

The description above applies to more general pairs of elements of $\mathfrak S$, namely whenever $S$ is replaced by some $Z$ so that $Y$ is properly \emph{nested into} $Z$, denoted $Y\propnest Z$. Nesting just means that the factor $Y$ is contained in the factor $Z$, or more precisely that there is a copy of $Y$ in a standard product region that is contained in a copy of $Z$.

\subsubsection{Orthogonality} We just saw that changing projections in both $Y$ and $Z$ when $Y\propnest Z$ can only be done in a rather specific way. The opposite situation is when $Y$ and $Z$ (which, recall, are factors of standard product regions) are \emph{orthogonal}, denoted $Y\orth Z$, meaning that they are (contained in) different factors of the same product region.  In this case, along a hierarchy path there is no constraint regarding which projection needs to change first, and in fact they can also change simultaneously. Ortohogonality is what creates non-hyperbolic behaviour in HHSs, and is what one has to constantly fight against.

\subsubsection{Transversality} When $Y,Z$ are neither $\nest$- nor $\orth$-comparable, we say that they are \emph{transverse}. This is the generic case. When $Y\transverse Z$ and $x,y\in\cuco X$ are so that $\dist_Y(x,y),\dist_Z(x,y)$ are both large, up to switching $Y,Z$ what happens is the following. When moving along any hierarchy path from $x$ to $y$ one has to first change the projection to $\fontact Y$ until it coarsely coincides with $\pi_Y(z)$, and only then the projection to $\fontact Z$ can start moving from $\pi_Z(x)$ to $\pi_Z(y)$.

Arguably the most useful feature of transversality is a slight generalisation of this. Given $x,y\in\cuco X$, and a set $\mathcal Y\subseteq \mathfrak S$ of pairwise transverse elements so that $\dist_{Y}(x,y)$ is large for every $Y\in\mathcal Y$, there is a total order on $\mathcal Y$ so that, whenever $Y<Z$, along any hierarchy path from $x$ to $y$ the projection to $Y$ has to change before the projection to $Z$ does, as described above.

\subsection{Realisation}

Even though we did not formally describe them, we saw that for certain pairs $Y,Z\in\mathfrak S$, namely when $Y\notorth Z$, there are some constraints on the projections of points in $\cuco X$ to $\fontact Y,\fontact Z$. These are called \emph{consistency inequalities}. As it turns out, the consistency inequalities are the only obstructions for ``coordinates'' $(b_Y\in\fontact Y)_{Y\in\mathfrak S}$ to be coarsely realised by a point $x$ in $\cuco X$, meaning that $\pi_Y(x)$ coarsely coincide with $b_Y$ in each $\fontact Y$. This is important because it allows to perform constructions in each of the $\fontact Y$ separately and then put everything back together. 

To make this principle clear, we now give an example of a construction of this type. Say we want to construct a ``coarse median'' map $m:\cuco X^3\to\cuco X$ (in the sense of \cite{Bow_coarse_median}), which let's just take to mean a coarsely Lipschitz map so that $m(x,x,y)$ is coarsely $x$. Consider $x,y,z$ in $\cuco X$, and let us define $m(x,y,z)$ by defining its coordinates in the $\fontact Y$. Given $Y\in\mathfrak S$, the triangle with vertices $\pi_Y(x),\pi_Y(y),\pi_Y(z)$ has a coarse centre $b_Y$, because $\fontact Y$ is hyperbolic. It turns out that the coordinates $(b_Y)$ satisfy the consistency inequalities, so that one can define $m(x,y,z)$ as the realisation point. As an aside, it is a nice exercise to use the properties of $m$ to show that $\cuco X$ satisfies a quadratic isoperimetric inequality.

To sum up, the distance formula says that the natural map $\cuco X\to\prod_{Y\in\mathfrak S}\fontact Y$ is ``coarsely injective'', and the consistency inequalities provide a coarse characterisation of the image.

\subsection{In the examples}

\subsubsection{RAAGs} 

In the case of RAAGs, $\mathfrak S$ (the set of factors of product regions) is the set of cosets of sub-RAAGs, considered up to \emph{parallelism}. We say that $gA_\Lambda, hA_\Lambda \subseteq A_{\Gamma}$ are parallel if $g^{-1}h$ commutes with every element of $A_\Lambda$, which essentially means that there's a product $g(A_\Lambda\times <g^{-1}h>)$ inside the RAAG $A_\Lambda$ so that $gA_\Lambda, hA_\Lambda$ are copies of one of the factors. Taking parallelism classes ensures that we will not do multiple counting in the distance formula. What we mean is that infinitely many parallel cosets would give the same contribution to the distance formula, which would clearly break it.

As in the case of $\fontact S$, $\fontact (gA_\Lambda)$ is a copy of the Cayley graph of $A_\Lambda$ with respect to the generating set $V\Lambda\cup \{A_{\Lambda'}<A_\Lambda:\Lambda'{\rm\ proper\ subgraph\ of\ }\Lambda\}$. The projection map from $A_\Gamma$ to $\fontact (gA_{\Gamma})$ is the composition of the closest-point projection to $gA_\Lambda$ in the usual Cayley graph of $A_\Lambda$, and the inclusion $gA_\Lambda\subseteq \fontact (gA_{\Gamma})$. The closest-point projection can also be rephrased in terms of the normal form for elements of $A_\Gamma$, since the normal form gives geodesics.

Nesting is inclusion up to parallelism, meaning that we declare $[gA_\Lambda]\nest [gA_{\Lambda'}]$ when $\Lambda\subseteq \Lambda'$, where $[\cdot]$ denotes the parallelism class. Similarly, we declare $[gA_\Lambda]\orth [gA_{\Lambda'}]$ if $\Lambda,\Lambda'$ form a join.

In the case of RAAGs, it turns out that geodesics in (the usual Cayley graph of $A_\Lambda$) are actually hierarchy paths.

\subsubsection{Mapping class groups} In this case, $\mathfrak S$ is the collection of (isotopy classes of essential) subsurfaces, with each $\fontact Y$ being the corresponding curve complex, and the maps $\pi_Y$ are defined using the so-called subsurface projections.
 Nesting is containment (up to isotopy), while orthogonality corresponds to disjointness (again up to isotopy).
 
 \subsubsection{CAT(0) cube complexes} Consider a CAT(0) cube complex $\cuco X$ with a factor system $\mathcal F$. In this case, $\mathfrak S$ is the union of $\{S=\cuco X\}$ and the set of \emph{parallelism classes} in $\mathcal F$. Parallelism can be defined in at least two equivalent ways. The first one is that the convex subcomplexes $F,F'$ are parallel if they cross the same hyperplanes. The second one, which provides a much better picture, is that $F,F'$ are parallel if there exists an isometric embedding of $F\times [0,n]\to \cuco X$, where $[0,n]$ is cubulated by unit intervals and $F\times [0,n]$ is regarded as a cube complex, so that $F\times\{0\}$ maps to $F$ in the obvious way, and the image of $F\times \{n\}$ is $F'$. As in the case of RAAGs, if we did not take parallelism classes then the distance formula would certainly not work due to multiple counting.
 
 The $\fontact [F]$ are obtained starting from $F$ and coning off all $F'\in\mathcal F$ contained in $F$. The maps $\pi_{[F]}$ are defined using gates.
 
 Nesting $[F]\nest [F']$ is inclusion up to parallelism, which can also be rephrased as: all hyperplanes crossing $F$ also cross $F'$ (notice that this does not depend on the choice of representatives). Orthogonality $[F]\orth [F']$ means that, up to parallelism, $F\times F'$ has a natural embedding into $\cuco X$. It can also be rephrased as: each hyperplane crossing $F$ crosses each hyperplane crossing $F'$.

%\section{Some constructions}

\part{Technical discussion}\label{part:technical}

Keeping in mind the heuristic discussion from Part \ref{part:heuristic}, we now analyse in more detail the definition and the main tools to study HHSs. We start with the axioms.

We will often motivate the axioms in terms of standard product regions, but we warn the reader in advance that those will be constructed only after we discuss all the axioms and a few tools. This, however, is inevitable. In fact, we are trying to describe a space that has some sort of subspaces, the standard product regions, that can be endowed with the same structure as the space itself. Until we know what that structure is in detail, we cannot use it to construct the standard product regions starting from first principles. Hopefully, one or more of the examples we discussed in Part \ref{part:heuristic} can help with intuition.

\section{Commentary on the axioms}
 
We will work in the context of a \emph{quasigeodesic space},  $\cuco X$, 
i.e., a metric space where any two
points can be connected by a uniform-quality quasigeodesic. It is more convenient for us to work with quasi-geodesic metric spaces than geodesic metric space because the standard product regions are in a natural way quasi-geodesic metric spaces, rather than geodesic metric space. Any quasi-geodesic metric space is quasi-isometric to a geodesic metric space since one can consider an approximating graph whose vertices form a maximal net, so for the purposes of large-scale geometry there's basically no difference between geodesic and quasi-geodesic metric spaces.

Actually, all the requirements in the definition of HHS are meant to be stable under passing to standard product regions. We do not have standard product regions yet, so what happens in the definition instead is that the axioms are about certain sub-collections of the set of hyperbolic spaces involved in the HHS structure, rather than just the whole collection.

We now go through the definition of HHS given in \cite{HHS2}, which is the one with ``optimised'' axioms compare to \cite{HHS1}. The statements of the axioms are given exactly as in \cite{HHS2}.

\par\bigskip

\emph{The $q$--quasigeodesic space  $(\cuco X,\dist_{\cuco X})$ is a \emph{hierarchically hyperbolic space} if there exists $\delta\geq0$, an index set $\mathfrak S$, and a set $\{\fontact W:W\in\mathfrak S\}$ of $\delta$--hyperbolic spaces $(\fontact U,\dist_U)$,  such that the following conditions are satisfied:
\begin{enumerate}
\item\textbf{(Projections.)}\label{item:dfs_curve_complexes} There is a set $\{\pi_W:\cuco X\rightarrow2^{\fontact W}\mid W\in\mathfrak S\}$ of \emph{projections} sending points in $\cuco X$ to sets of diameter bounded by some $\xi\geq0$ in the various $\fontact W\in\mathfrak S$.  Moreover, there exists $K$ so that each $\pi_W$ is $(K,K)$--coarsely Lipschitz.\end{enumerate}}

\par\bigskip

The index set $\mathfrak S$ is the set of factors of standard product regions. Any $V\in \mathfrak S$ hence corresponds to each of many ``parallel'' subsets of $\cuco X$.
We already saw where the hyperbolic spaces associated to an HHS comes from: each factor of a standard product region contains various standard product sub-regions, which we can cone-off to obtain a hyperbolic space. 
The way to think about the projection is that the standard product regions and their factors come with a coarse retraction from $\cuco X$, and the projections $\pi_W$ in the definition are the composition of those retractions with the cone-off map. This is admittedly a bit circular because we will later define the retractions in terms of the $\pi_W$, but should hopefully help to understand the picture.

The reason why the projections take value in bounded subsets of the $\fontact W$ rather than points is just that in several situations, for example subsurface projections for mapping class groups, this is what one gets in a natural way. One can make arbitrary choices and modify the projections to take value in points, and nothing would be affected.

\par\bigskip

\emph{\begin{enumerate}[(2)] \item \textbf{(Nesting.)} \label{item:dfs_nesting} $\mathfrak S$ is
 equipped with a partial order $\nest$, and either $\mathfrak
 S=\emptyset$ or $\mathfrak S$ contains a unique $\nest$--maximal
 element; when $V\nest W$, we say $V$ is \emph{nested} in $W$.  We
 require that $W\nest W$ for all $W\in\mathfrak S$.  For each
 $W\in\mathfrak S$, we denote by $\mathfrak S_W$ the set of
 $V\in\mathfrak S$ such that $V\nest W$.  Moreover, for all $V,W\in\mathfrak S$
 with $V\propnest W$ there is a specified subset
 $\rho^V_W\subset\fontact W$ with $\diam_{\fontact W}(\rho^V_W)\leq\xi$.
 There is also a \emph{projection} $\rho^W_V\colon \fontact
 W\rightarrow 2^{\fontact V}$.  (The similarity in notation is
 justified by viewing $\rho^V_W$ as a coarsely constant map $\fontact
 V\rightarrow 2^{\fontact W}$.)\end{enumerate}}
 
\par\bigskip
 
 Nesting corresponds to inclusion between standard product regions. The maximal element corresponds to $\cuco X$ itself, thought of as a product region with a trivial factor.
 
 Recall that the $\fontact W$ are obtained coning-off standard product regions, i.e. making them bounded. For $V\propnest W$, the bounded set $\rho^V_W$ is one such bounded set, where $V$ is regarded as a standard product region with one trivial factor. In the other direction, $\rho^W_V$ is obtained restricting the retraction to $W$.
 
 We will discuss below the fact that $\rho^W_V$ for $V\propnest W$ is not strictly needed, and can be to all effects and purposes be reconstructed from $\pi_W$ and $\pi_V$.
 
 Regarding the notation, the $\rho^V_W$s in this axiom as well as the ones below always go ``from top to bottom'', meaning that $\rho^V_W$ is always some kind of projection from $V$ to $W$.
 
 \par\bigskip

  \emph{\begin{enumerate}[(3)]
 \item \textbf{(Orthogonality.)} 
 \label{item:dfs_orthogonal} $\mathfrak S$ has a symmetric and anti-reflexive relation called \emph{orthogonality}: we write $V\orth W$ when $V,W$ are orthogonal.  Also, whenever $V\nest W$ and $W\orth U$, we require that $V\orth U$.  Finally, we require that for each $T\in\mathfrak S$ and each $U\in\mathfrak S_T$ for which $\{V\in\mathfrak S_T:V\orth U\}\neq\emptyset$, there exists $W\in \mathfrak S_T-\{T\}$, so that whenever $V\orth U$ and $V\nest T$, we have $V\nest W$.  Finally, if $V\orth W$, then $V,W$ are not $\nest$--comparable.
 \end{enumerate}}

 \par\bigskip
 
 Orthogonality is what creates non-trivial products: $V$ and $W$ are orthogonal if they participate in a common standard product region, meaning that they are distinct factors. With this interpretation, it should be clear why when $V\nest W$ and $W\orth U$, we require $V\orth U$, and also why when $V\orth W$, then $V,W$ should not be $\nest$--comparable.
 
 The tricky part is the one about $\{V\in\mathfrak S_T:V\orth U\}$. Let us first discuss the case of $\mathfrak S$ instead of more general $\mathfrak S_T$. The point is that one wants to define an orthogonal complement of the $V\in\mathfrak S$, and one wants it to be an HHS, with corresponding index set $U^\perp=\{U\in\mathfrak S:V\orth U\}$. For that to be the case, one would want $U^\perp$ to contain a $\nest$--maximal element (if it is non-empty). The axiom says something a bit weaker, because the $W$ containing each $V\orth U$ is not required to be itself orthogonal to $U$. This is still enough to have an HHS structure on the orthogonal complement. The only reason we did not require the stronger version with $W\perp U$ in \cite{HHS1} is that at the time we were not able to prove that such $W$ exists in the case of CAT(0) cube complexes. However, it follows from \cite[Theorem 3.5]{HagenSusse} that proper cocompact cube complexes satisfy the stronger version of the axiom, so that in fact all natural examples of HHS (so far) do, and there is no harm in strengthening the orthogonality axiom. In fact, sometimes the weaker formulation gives technical problems.
 
 As a final comment, it is natural to formulate the axiom for general $\mathfrak S_T$ instead of just for $\mathfrak S$ because all axioms need to work inductively for product regions.
 
\par\bigskip
 
 {\it \begin{enumerate}[(4)]
 \item \textbf{(Transversality and consistency.)}
 \label{item:dfs_transversal} If $V,W\in\mathfrak S$ are not
 orthogonal and neither is nested in the other, then we say $V,W$ are
 \emph{transverse}, denoted $V\transverse W$.  There exists
 $\kappa_0\geq 0$ such that if $V\transverse W$, then there are
  sets $\rho^V_W\subseteq\fontact W$ and
 $\rho^W_V\subseteq\fontact V$ each of diameter at most $\xi$ and 
 satisfying: $$\min\left\{\dist_{
 W}(\pi_W(x),\rho^V_W),\dist_{
 V}(\pi_V(x),\rho^W_V)\right\}\leq\kappa_0$$ for all $x\in\cuco X$.
 
 For $V,W\in\mathfrak S$ satisfying $V\nest W$ and for all
 $x\in\cuco X$, we have: $$\min\left\{\dist_{
 W}(\pi_W(x),\rho^V_W),\diam_{\fontact
 V}(\pi_V(x)\cup\rho^W_V(\pi_W(x)))\right\}\leq\kappa_0.$$ 
 
 The preceding two inequalities are the \emph{consistency inequalities} for points in $\cuco X$.
 
 Finally, if $U\nest V$, then $\dist_W(\rho^U_W,\rho^V_W)\leq\kappa_0$ whenever $W\in\mathfrak S$ satisfies either $V\propnest W$ or $V\transverse W$ and $W\not\orth U$.
 \end{enumerate}}

 \par\bigskip

 Transversality is best thought of as being in ``general position''. As an aside for the reader who speaks relative hyperbolicity, if $U\transverse V$, then they behave very similarly to distinct cosets of peripheral subgroups of a relatively hyperbolic group; for example the  projections to $U,V$ should be compared to closest-point projections onto a pair of distinct cosets. The first consistency inequality, also known as Behrstock inequality, is very important, so we now discuss a few ways to think about it (and its consequences). Incidentally, we note that the Behrstock inequality is important beyond the HHS world too; for example, it plays a prominent role in the context of the projection complexes from \cite{BBF}, which have many applications.
 
 In words, the Behrstock inequality says that if $V\transverse W$ and $x\in\cuco X$ projects far from $\rho^W_V$ in $\fontact V$, then $x$ projects close to $\rho^V_W$ in $\fontact W$. ($\rho^W_V$ is best thought of as the projection of $W$ onto $\fontact V$.)
 Let us start by discussing an easy situation where the inequality holds. Suppose that $V$ and $W$ are two quasi-convex subsets of a hyperbolic space and suppose that $\pi_V(W),\pi_W(V)$ are both bounded, where $\pi_V,\pi_W$ denote (coarse) closest point projections. Then setting $\rho^V_W=\pi_V(W)$ and $\rho^W_V=\pi_W(V)$, the first consistency inequality holds, and it is illustrated in the following picture:
 
 \begin{figure}[h]
  \includegraphics[scale=0.8]{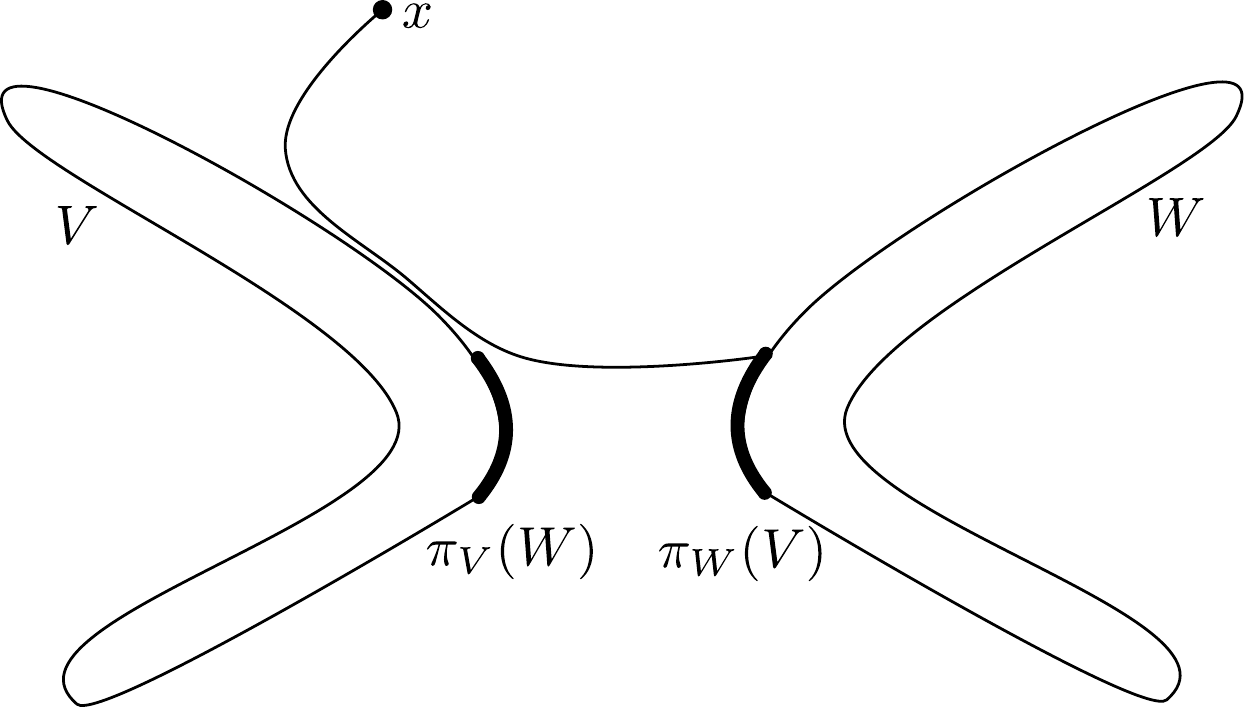}
  \caption{The Behrstock inequality for quasiconvex subspaces of a hyperbolic space.}\label{fig:Behr}
 \end{figure}

 Here is a sketch of the argument, which should also clarify the meaning of the inequality. If $x$ projects to $V$ far away from the projection of $W$, as in the picture, the we have to show that it projects close to the projection to $V$ onto $W$. This is because any geodesic from $x$ to $W$ must pass close to $V$, by a standard hyperbolicity argument.
 
 This last fact is useful to keep in mind: in the situation above, to go from $x$ to $W$ one has to pass close to $V$ first, and change the projection to $V$ in the process.
 
 A second way to understand the inequality is to draw the image of $\pi=\pi_V\times\pi_W$, which is coarsely the following ``cross'':
 
  \begin{figure}[h]
  \includegraphics[scale=0.8]{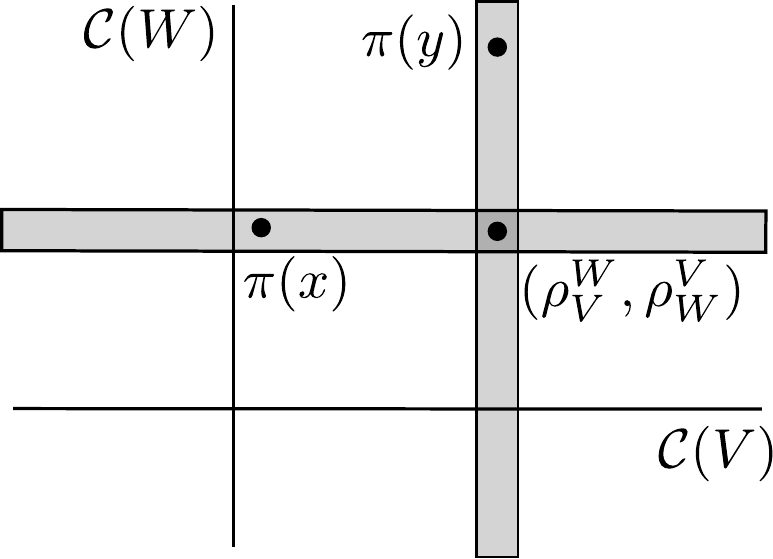}
  \caption{The image of $\pi=\pi_V\times\pi_W$, when $V\transverse W$.}\label{fig:cross}
 \end{figure}
 
 From this graph we see a similar phenomenon to the one above: depending on where $\pi(x),\pi(y)$ lie on the cross, to go from $x$ to $y$ one has to change the projection to $V$ first or the projection to $W$ first.
 
 This brings us to an important consequence of the Behrstock inequality, which is that one can order transverse $V,W$ that lie ``between'' $x$ and $y$. Suppose that $x,y\in \cuco X$ and $\{V_i\}$ are pairwise transverse and so that $\dist_{V_i}(x,y)$ are all much larger than the constant in the Behrstock inequality. Then for each $V,W\in\{V_i\}$, up to switching $V,W$, the situation looks like Figure \ref{fig:cross}, and in this case we write $V\prec W$. We can give several equivalent description of ``the picture above'', and manipulating the Behrstock inequality reveals that they are all equivalent. These are the following, where we are assuming $\dist_V(\pi_V(x),\pi_V(y)),\dist_W(\pi_W(x),\pi_W(y))\geq 10E$ for some sufficiently large $E$:
 \begin{itemize}
  \item $V\prec W$,
  \item $\dist_W(\pi_W(x),\rho^V_W)\leq E$,
  \item $\dist_W(\pi_W(y),\rho^V_W)> E$,
  \item $\dist_V(\pi_V(x),\rho^W_V)> E$,
  \item $\dist_V(\pi_V(y),\rho^W_V)\leq E$.
 \end{itemize}
A very important fact is that $\prec$ is a total order on $\{V_i\}$. My favourite way to draw this is the following, assuming for simplicity $V_i\prec V_j$ if and only if $i<j$:

  \begin{figure}[h]
  \includegraphics[width=\textwidth]{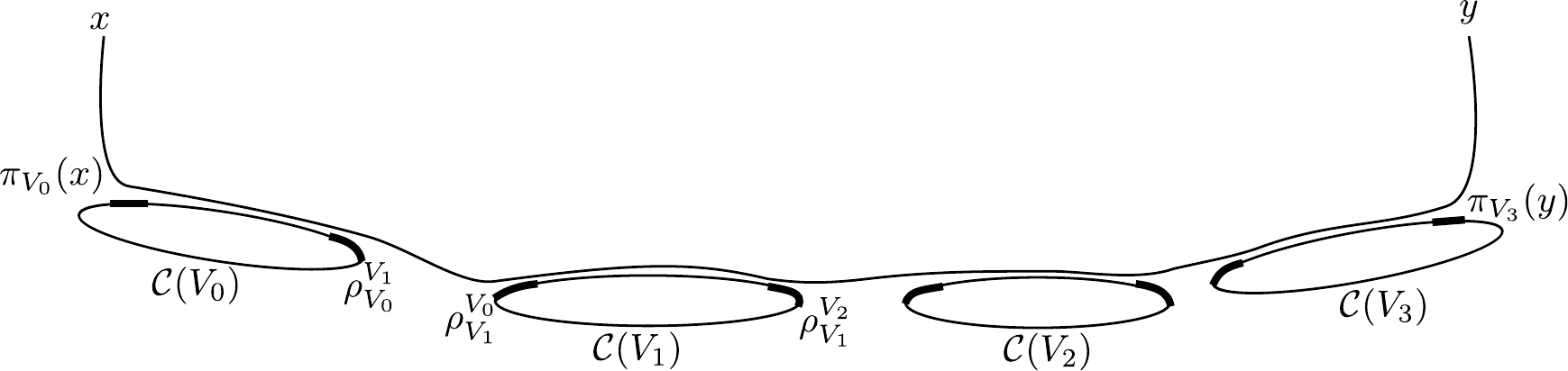}
 \end{figure}

This picture does not really take place anywhere, but it contains interesting information. You can pretend that the $\fontact{V_i}$ are quasiconvex subsets of a hyperbolic space as in Figure \ref{fig:Behr}, with the path from $x$ to $y$ in the picture representing a geodesic from $x$ to $y$ that passes close to them in the order given by $\prec$. From the picture you can read off where the various $\rho$s are by following the path. In particular, you see that for $i<j<k$, $\rho^{V_j}_{V_i}$ and $\rho^{V_k}_{V_i}$ coarsely coincide with each other and with $\pi_V(y)$. This picture still works to understand where projections lie if you, for example, add another point $z$. You can try to convince yourself, first from the picture and then formally, that if $z$ projects ``in the middle'' on some $\fontact V_i$ then, for $j>i$, $\pi_{V_j}(z)$ coarsely coincides with $\pi_{V_j}(x)$.

We now discuss the second consistency inequality in conjunction with another axiom:
\par\smallskip
 \emph{\begin{enumerate}[(7)]
 \item \textbf{(Bounded geodesic image.)} \label{item:dfs:bounded_geodesic_image} For all $W\in\mathfrak S$, all $V\in\mathfrak S_W-\{W\}$, and all geodesics $\gamma$ of $\fontact W$, either $\diam_{\fontact V}(\rho^W_V(\gamma))\leq E$ or $\gamma\cap\neb_E(\rho^V_W)\neq\emptyset$. 
 \end{enumerate}}

 \par\bigskip
 In words, when $V$ is properly nested into $W$, then the projection $\rho^W_V$ from $\fontact W$ to $\fontact V$ is coarsely constant along geodesics far from $\rho^V_W$ (recall that this is the copy of $V$ that gets coned-off to make $\fontact W$ out of $W$).
 
 This is virtually always used together with the second consistency inequality, which implies that if $\pi_W(x)$ is far from $\rho^V_W$ for some $x\in\cuco X$ then $\rho^W_V(\pi_W(x))$ coarsely coincides with $\pi_V(x)$. This yields the version of bounded geodesic image that most often gets used in practice:
 
 \begin{lem}[{See e.g. \cite[Lemma 1.5]{HHS_quasiflats}}]
  Let $(X,\mathfrak S)$ be hierarchically hyperbolic. Up to increasing $E$ as in the bounded geodesic image axiom, for all $W\in\mathfrak S$, all $V\in\mathfrak S_W-\{W\}$, and all $x,y\in\cuco X$ so that some geodesic from $\pi_W(x)$ to $\pi_W(y)$ stays $E$--far from $\rho^V_W$, we have $\dist_V(\pi_V(x),\pi_V(y))\leq E$.
 \end{lem}
 
  \begin{figure}[h]
  \includegraphics[scale=0.7]{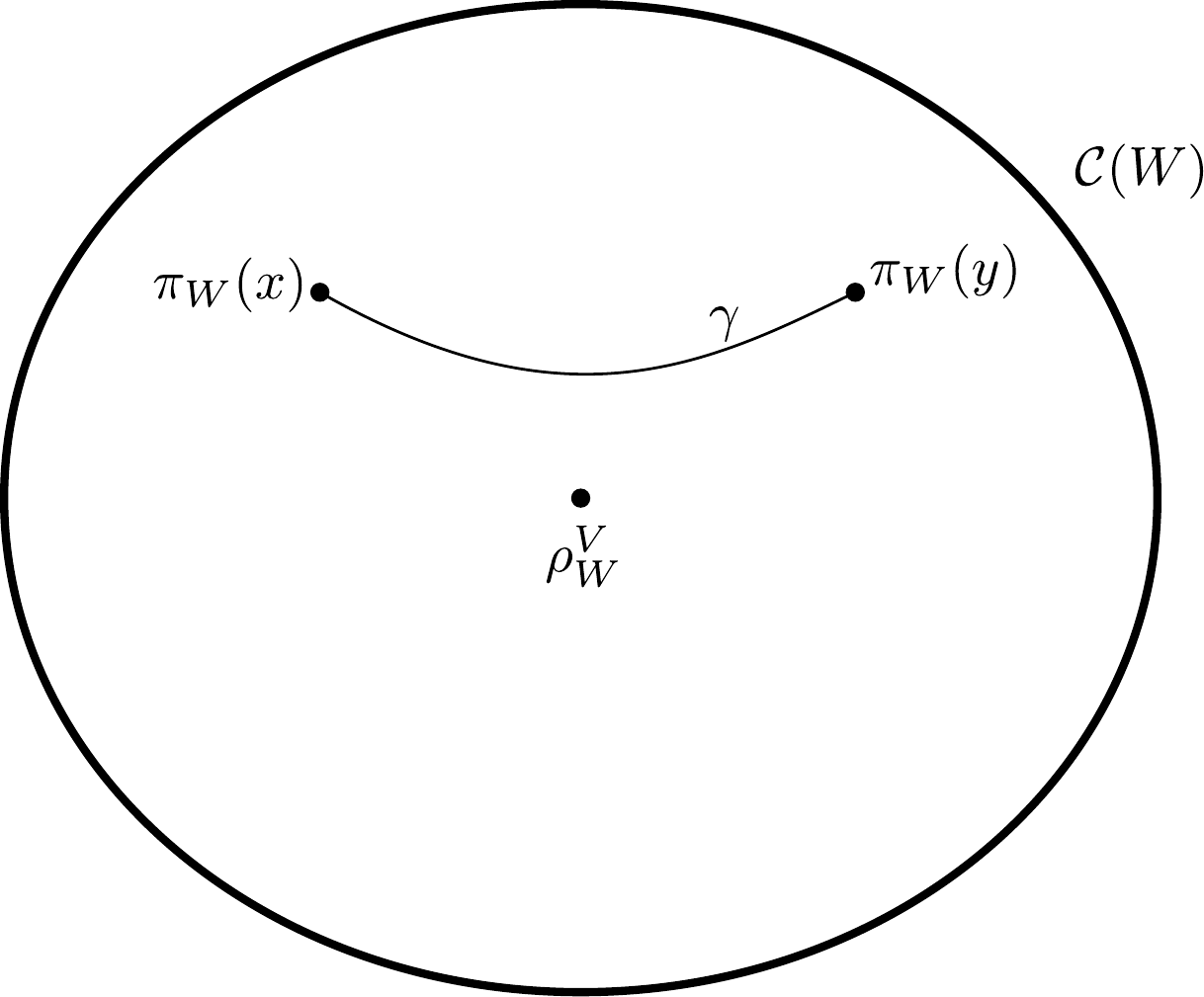}
  \caption{In the picture we have $V\propnest W$ and $\gamma$ is a geodesic. According to bounded geodesic image, $\pi_V(x)$ and $\pi_V(y)$ coarsely coincide.}
 \end{figure}

One can simply replace the bounded geodesic image axiom and the second consistency inequality with the lemma, since $\rho^W_V$ can be reconstructed from $\pi_W$ and $\pi_V$ at least on $\pi_W(\cuco X)$ in view of the lemma. However, for some purposes one still needs $\rho^W_V$. This is most notably the case for the realisation theorem.

Another picture to keep in mind regarding bounded geodesic image is that, given $x,y\in\cuco X$ and $W$, one can consider all $V\propnest W$ so that $\dist_V(\pi_V(x),\pi_V(y))$ is large. The corresponding $\rho^V_W$ will form a ``halo'' around a geodesic from $\pi_W(x)$ to $\pi_W(y)$.

  \par\bigskip
 
 \emph{\begin{enumerate}[(5)]
 \item \textbf{(Finite complexity.)} \label{item:dfs_complexity} There exists $n\geq0$, the \emph{complexity} of $\cuco X$ (with respect to $\mathfrak S$), so that any set of pairwise--$\nest$--comparable elements has cardinality at most $n$.
 \end{enumerate}}
 
 \par\bigskip
 
 This axiom should be pretty self-explanatory. Induction on complexity is very common in the HHS world. The base case (complexity $1$) is that of hyperbolic spaces.

\par\bigskip
 
 \emph{\begin{enumerate}[(6)]
 \item \textbf{(Large links.)} \label{item:dfs_large_link_lemma} There
 exists $E\geq\max\{\xi,\kappa_0\}$ such that the following holds.
 Let $W\in\mathfrak S$ and let $x,x'\in\cuco X$.  Let
 $N=\dist_{W}(\pi_W(x),\pi_W(x'))$.  Then, 
%  either
%  $$\dist_T(\pi_T(x),\pi_T(x'))\leq E$$ for all $T\in\mathfrak
%  S_W-\{W\}$, or 
 there exist $T_1,\ldots,T_{\lfloor
 N\rfloor}\in\mathfrak S_W-\{W\}$ such that for all $T\in\mathfrak
 S_W-\{W\}$, either $T\in\mathfrak S_{T_i}$ for some $i$, or $\dist_{
 T}(\pi_T(x),\pi_T(x'))<E$.  Also, $\dist_{
 W}(\pi_W(x),\rho^{T_i}_W)\leq N$ for each $i$.
 \end{enumerate}}

\par\bigskip

In words, the axioms say that, given $W$ and $x,x'\in\cuco X$, each of the $V\propnest W$ so that $\dist_V(\pi_V(x),\pi_V(x'))$ is nested into one of a few fixed $T_i\propnest W$. The number of $T_i$ required is bounded only in terms of $\dist_W(\pi_W(x),\pi_W(x'))$ (which can be much smaller than their distance in $\cuco X$).

This axiom is very related to bounded geodesic image, and in fact in concrete examples they are often proven at the same time. Bounded geodesic image provides a ``halo'' of $\rho^V_W$ around a geodesic connecting $\pi_W(x),\pi_W(x')$, and there can be arbitrarily many of these. However, large links organises them into a few (possibly intersecting) subsets, each of which contains the $\rho^V_W$ with $V$ nested into some fixed $T_i$. The number of such $T_i$ is bounded in terms of the distance $\dist_W(\pi_W(x),\pi_W(x'))$.

Large links is used in arguments of the following type. Consider two points $x,y$ that are far in $\cuco X$. If they are far in $\fontact S$ (meaning that their projections are), then one can use the geometry of $\fontact S$ to study whatever property one is interested in. Otherwise, there are few $T_i$s, and one can then analyse corresponding standard product regions. In one of them, (the retractions of) $x,y$ are still far away, so one can use induction based on the fact that the standard product region is an HHS of strictly lower complexity.

One concrete lemma that makes this more precise is the "passing up" lemma \cite[Lemma 2.5]{HHS2}. This says the (contrapositive of the) following. If one has $x,y\in\cuco X$ and some $S_i\in\mathfrak S$ so that the $\dist_{S_i}(\pi_{S_i}(x), \pi_{S_i}(y))$ are all large and each $S_i$ is $\nest$--maximal with this property, then there is a bound on how many $S_i$ there are.
 
  \par\bigskip
 
 \emph{\begin{enumerate}[(8)]
 \item \textbf{(Partial Realization.)} \label{item:dfs_partial_realization} There exists a constant $\alpha$ with the following property. Let $\{V_j\}$ be a family of pairwise orthogonal elements of $\mathfrak S$, and let $p_j\in \pi_{V_j}(\cuco X)\subseteq \fontact V_j$. Then there exists $x\in \cuco X$ so that:
 \begin{itemize}
 \item $\dist_{V_j}(x,p_j)\leq \alpha$ for all $j$,
 \item for each $j$ and 
 each $V\in\mathfrak S$ with $V_j\propnest V$, we have 
 $\dist_{V}(x,\rho^{V_j}_V)\leq\alpha$, and
 \item if $W\transverse V_j$ for some $j$, then $\dist_W(x,\rho^{V_j}_W)\leq\alpha$.
 \end{itemize}
 \end{enumerate}}

 \par\bigskip
 
 Roughly speaking, the axiom says that, given pairwise-orthogonal $\{V_i\}$ there is no restriction on the projections of points of $\cuco X$ to the $\fontact V_i$; any choice of coordinates can be realised by a point in $\cuco X$. This is the opposite of what happens when $V\not\orth W$ (i.e. in one of the cases $V=W,V\propnest W,W\propnest V$ or $V\transverse W$), in which case there are serious restrictions on the projections in view of the consistency inequalities. 
 
 This axiom gives us the first glimpse of how the standard product regions arise, and what their coordinates in the various $\fontact U$ look like. Starting from the family of pairwise orthogonal elements  $\{V_j\}$, we see that the axioms provides us with the freedom to move independently in each of the $\fontact V_j$. When we will have the ``full'' realisation theorem, this will give us a product region associated to $\{V_j\}$. The second condition can be explained as follows: the coordinate in $\fontact V$ does not coarsely vary when moving around the standard product region because the standard product region is coned-off there. The third condition tells us that ``generic'' pairs of standard product regions do not interact much with each other. (Recall that we think of transversality as being in ``generic position''.)  

  \par\bigskip
 
 \emph{\begin{enumerate}[(9)]
\item\textbf{(Uniqueness.)} For each $\kappa\geq 0$, there exists
$\theta_u=\theta_u(\kappa)$ such that if $x,y\in\cuco X$ and
$\dist(x,y)\geq\theta_u$, then there exists $V\in\mathfrak S$ such
that $\dist_V(x,y)\geq \kappa$.\label{item:dfs_uniqueness}
\end{enumerate}}

\par\bigskip

Informally, the axiom says that if $x,y$ are close in each $\fontact V$ (meaning that their projections are) then $x,y$ are close in $\cuco X$.

This axiom is a weaker form of distance formula. The point is that it is in many circumstances much easier to prove than the ``full'' distance formula, allowing for easier proofs that certain spaces are HHS. This is the case for mapping class groups, where there's a one-page argument for this axiom, given in \cite[Section 11]{HHS2}, while the known proofs of the distance formula are much more involved.

\section{Main tools}

In addition to the axioms, there are 3 fundamental properties of HHSs. These were actually part of the first set of axioms, but they have a much higher level of sophistication than any of the axioms.

\subsection{Distance formula}

We stated the distance formula in Part \ref{part:heuristic}, but let us recall it. Given $A,B\in\reals$, the symbol $\ignore{A}{B}$ will denote $A$ if $A\geq B$ and $0$ otherwise.  Given $C,D$, we write $A\asymp_{C,D}B$ to mean $C^{-1}A-D\leq B\leq CA+D$.

To save notation, we denote $\dist_{ W}(x,y)=\dist_{ W}(\pi_W(x),\pi_W(y))$.

\begin{thm}[{Distance Formula, \cite[Theorem 4.5]{HHS2}}]\label{thm:distance_formula}
 Let $(X,\mathfrak S)$ be hierarchically hyperbolic. Then there exists $s_0$ such that for all $s\geq s_0$ there exist
 constants $K,C$ such that for all $x,y\in\cuco X$,
 $$\dist_{\cuco X}(x,y)\asymp_{K,C}\sum_{W\in\mathfrak S}\ignore{\dist_{ W}(x,y)}{s}.$$
\end{thm}

The distance formula allows one to reconstruct the geometry of $\cuco X$ from that of the hyperbolic spaces $\fontact W$, and at this point its importance should hopefully be evident. It is important to note that the distance formula works for any sufficiently high threshold. This is useful in practice because typically one proceeds along the following lines. One starts with a configuration in $\cuco X$, projects it to the $\fontact W$ and keeps into account the distance formula to figure out what one gets. Then one performs some coarse construction in the $\fontact W$, and then goes back to $\cuco X$. In the process, more often than not some projections gets moved a bounded amount. To compensate for this, one uses a higher threshold in the distance formula.

\subsection{Hierarchy paths}

%A \emph{$(D,D)$--quasigeodesic} in the metric space $M$ is a
%$(D,D)$--quasi-isometric embedding $f\co[0,\ell]\to M$; we allow $f$ to
%be a coarse map, i.e., to send points in $[0,\ell]$ to uniformly
%bounded sets in $M$.

Hierarchy paths are quasigeodesics in $\cuco X$ that shadow geodesics in all $\fontact W$, which is clearly a very nice property to have since we want to relate the geometry of $\cuco X$ to that of the $\fontact W$. Let us define them precisely.

For $M$ a metric space, a (coarse) map $f\co[0,\ell]\to M$ is a
\emph{$(D,D)$--unparameterized quasigeodesic} if there exists a
strictly increasing function $g\co[0,L]\to[0,\ell]$ such that $f\circ
g\co[0,L]\to M$ is a $(D,D)$--quasigeodesic and for each
$j\in[0,L]\cap\naturals$, we have $\diam_M\left(f(g(j))\cup
f(g(j+1))\right)\leq D$.

\begin{defn}[Hierarchy path]\label{defn:hierarchy_path}
Let $(X,\mathfrak S)$ be hierarchically hyperbolic. For $D\geq 1$, a (not necessarily continuous) path $\gamma\co[0,\ell]\to\cuco X$ is a \emph{$D$--hierarchy path} if
 \begin{enumerate}
  \item $\gamma$ is a $(D,D)$-quasigeodesic,
  \item for each $W\in\mathfrak S$, the path $\pi_W\circ\gamma$ is an unparameterized $(D,D)$--quasigeodesic.
\end{enumerate}
\end{defn}

\begin{thm}[{Existence of Hierarchy Paths, \cite[Theorem 4.4]{HHS2}}]\label{thm:monotone_hierarchy_paths}
Let $(\cuco X,\mathfrak S)$ be hierarchically hyperbolic. Then there exists $D_0$ so that any $x,y\in\cuco X$ are joined by a $D_0$-hierarchy path. 
\end{thm}

Whenever possible, one should work with hierarchy paths rather than other quasigeodesics, even actual geodesics. Unfortunately, not all quasigeodesics are hierarchy paths (meaning that one cannot control how close the projection to some $\fontact W$ of a $(D,D)$--quasigeodesic is to being a geodesic as a function of $D$ only). In fact, there are spiraling quasigeodesics in $\mathbb R^2$, and, even worse than that, it is a folklore result that in mapping class groups there are quasigodesics that project to ``arbitrarily bad'' paths even in the curve graph of the whole surface.

Moreover, hierarchy paths with given endpoints are not coarsely unique: think of $\mathbb R^2$, where there are plenty of quasigeodesics monotone in each factor that connect points far away along a diagonal. In fact, it is a very important problem to study to what extent one can make hierarchy paths canonical by adding more restrictions.

\subsection{Realisation}

In this subsection we discuss the realisation theorem, which says that the consistency inequalities characterise the tuples $(\pi_W(x))_{W\in\mathfrak S}$ for $x\in\cuco X$. We think of the $\pi_W(x)$ as the coordinates of $x$.

\begin{defn}[Consistent]\label{defn:consistent}
Fix $\kappa\geq0$ and let $\tup b\in\prod_{U\in\mathfrak S}2^{\fontact U}$ be a tuple such that for each $U\in\mathfrak S$, the coordinate $b_U$ is a subset of $\fontact U$ with $\diam_{\fontact U}(b_U)\leq\kappa$.  The tuple $\tup b$ is \emph{$\kappa$--consistent} if, whenever $V\transverse W$,
$$\min\left\{\dist_{ W}(b_W,\rho^V_W),\dist_{ V}(b_V,\rho^W_V)\right\}\leq\kappa$$
and whenever $V\nest W$, 
$$\min\left\{\dist_{ W}(b_W,\rho^V_W),\diam_{\fontact V}(b_V\cup\rho^W_V(b_W))\right\}\leq\kappa.$$
\end{defn}

(Notice that in the definition of consistent tuple we need the map $\rho^W_V$ for $V\propnest W$.)

\begin{thm}[{Realisation of consistent tuples, \cite[Theorem 3.1]{HHS2}}]\label{thm:realization}
 For each $\kappa\geq1$ there exist $\theta_e,\theta_u\geq0$ such that
 the following holds.  Let $\tup
 b\in\prod_{W\in\mathfrak S}2^{\fontact W}$ be $\kappa$--consistent;
 for each $W$, let $b_W$ denote the $\fontact W$--coordinate of $\tup
 b$.

 Then there exists $x\in \cuco X$ so that $\dist_{
 W}(b_W,\pi_W(x))\leq\theta_e$ for all $\fontact W\in\mathfrak S$.
 Moreover, $x$ is \emph{coarsely unique} in the sense that the set of
 all $x$ which satisfy $\dist_{ W}(b_W,\pi_W(x))\leq\theta_e$ in each
 $\fontact W\in\mathfrak S$, has diameter at most $\theta_u$.
\end{thm}

As mentioned in Part \ref{part:heuristic}, the realisation theorem is used to perform constructions in all the $\fontact W$ separately and then pull those back to $\cuco X$. One such construction is (at last!) that of standard product regions. Basically, we fix $U\in \mathfrak S$, and consider partial systems of coordinates $(b_V)$, where we only assign $b_V$ when either $V\nest U$ or $V\orth U$. If this partial system of coordinates satisfies the consistency inequalities, we can extend it and use realisation to find a corresponding point in $\cuco X$. The standard product region associated to $U$ is obtained considering all such realisation points. A similar game can be played starting from pairwise orthogonal $U_i$, but for simplicity we stick to the case of a single $U$. Let us make this more precise.

\begin{defn}[Nested partial tuple]\label{defn:nested_partial_tuple}
Recall $\mathfrak S_U=\{V\in\mathfrak S:V\nest U\}$.  Fix
$\kappa\geq\kappa_0$ and let $\mathbf F_U$ be the set of
$\kappa$--consistent tuples in $\prod_{V\in\mathfrak S_U}2^{\fontact
V}$.
\end{defn}

\begin{defn}[Orthogonal partial tuple]\label{defn:orthogonal_partial_tuple}
Let $\mathfrak S_U^\orth=\{V\in\mathfrak S:V\orth U\}\cup\{A\}$, where
$A$ is a $\nest$--minimal element $W$ such that $V\nest W$ for all
$V\orth U$.  Fix $\kappa\geq\kappa_0$, let $\mathbf E_U$ be the set of
$\kappa$--consistent tuples in $\prod_{V\in\mathfrak
S_U^\orth-\{A\}}2^{\fontact V}$.
\end{defn}

\begin{cons}[Product regions in $\cuco X$]\label{const:embedding_product_regions}
Given $\cuco X$ and $U\in\mathfrak S$, there are coarsely well-defined
maps $\phi^\nest,\phi^\orth\co\mathbf F_U,\mathbf E_U\to\cuco X$ which extend to a coarsely
well-defined map $\phi_U\co\mathbf F_U\times \mathbf E_U\to\cuco X$, whose image $P_U$ we call a \emph{standard product region}. Indeed, for each $(\tup a,\tup b)\in
\mathbf F_U\times \mathbf E_U$, and each $V\in\mathfrak S$, define the
co-ordinate $(\phi_U(\tup a,\tup b))_V$ as follows.  If $V\nest U$,
then $(\phi_U(\tup a,\tup b))_V=a_V$.  If $V\orth U$, then
$(\phi_U(\tup a,\tup b))_V=b_V$.  If $V\transverse U$, then
$(\phi_U(\tup a,\tup b))_V=\rho^U_V$.  Finally, if $U\nest V$, and
$U\neq V$, let $(\phi_U(\tup a,\tup b))_V=\rho^U_V$.

By design of the axioms, it is straightforward (but a bit tedious) to check that we actually defined a consistent tuple, see \cite[Section 13.1]{HHS1}.
\end{cons}

We notice that by the very definition of $P_U$, the following hold. First, $\pi_Y(P_U)$ is uniformly bounded if $U\propnest Y$ (making sure that it makes sense to think of $U$ as being coned-off to get $\fontact Y$), as well as if $U\transverse Y$ (so that we can actually think of $P_U$ and $P_Y$ as ``independent'').

Coarse retractions onto standard product regions have been mentioned above. It should not be hard to guess how they are constructed at this point. One simply starts with $x\in \cuco X$, define coordinates by taking $\pi_Y(x)$ whenever $Y\nest U$ or $Y\orth U$ and $\rho^U_Y$ otherwise, and takes a realisation point. Basically, one defines the retraction of $x\in\cuco X$ by keeping the coordinates involved in the standard product region only.

\subsubsection{Coarse median}
As mentioned in Part \ref{part:heuristic}, the realisation theorem can be used to construct a coarse-median map in the sense of \cite{Bow_coarse_median} (also called centroid in \cite{BM_centroid}). This is the map $\median\co\cuco X^3\rightarrow\cuco X$ defined as follows.  Let $x,y,z\in\cuco X$ and, for each $U\in\mathfrak S$, let $b_U$ be a coarse centre for the triangle with vertices $\pi_U(x),\pi_U(y),\pi_U(z)$. More precisely, $b_U$ is any point in $\fontact U$ with the property that there exists
 a geodesic triangle in $\fontact U$ with vertices in
 $\pi_U(x),\pi_U(y),\pi_U(z)$ each of whose sides contains a point
 within distance $\delta$ of $b_U$, where $\delta$ is the hyperbolicity constant of $\fontact U$.

By \cite[Lemma 2.6]{HHS2} (which is easy, and a good exercise), the tuple $\tup b\in\prod_{U\in\mathfrak S}2^{\fontact U}$ whose $U$-coordinate is $b_U$ is $\kappa$--consistent for an appropriate choice of $\kappa$.  Hence, by the realisation theorem, there exists $\median=\median(x,y,z)\in\cuco X$ such that $\dist_U(\median,b_U)$ is uniformly bounded for all $U\in\mathfrak S$.  Moreover, this is coarsely well-defined, by the uniqueness axiom. The fact that this coarse median map actually makes the HHS into a coarse median space, and that, moreover, the rank is the ``expected'' one, is \cite[Corollary 2.15]{HHS_quasiflats}. 

The existence of the coarse median has many useful consequences, for example regarding asymptotic cones \cite{Bow_coarse_median,Bowditch:rigid} (\cite{HHS_quasiflats} heavily relies on these). Also, the explicit construction itself is useful in various arguments, for example to construct the kind of retractions mentioned in the subsection on hierarchical quasiconvexity below.

\section{Additional tools}

\subsection{Hierarchical quasiconvexity}

Quasiconvex subspaces are important in the study of hyperbolic spaces. The corresponding notion in the HHS world is \emph{hierarchical quasiconvexity}. Prominent examples of hierarchically quasiconvex subspaces are standard product regions.

The natural first guess for what a hierarchically quasiconvex subspace of an HHS $\cuco X$ should be is a subspace that projects to uniformly quasiconvex subspaces in all $\fontact U$. This is part of the definition, but not quite enough to have a good notion. In fact, the aforementioned property is satisfied by subspaces that are not even coarsely connected. There are at least two ideas to ``complete'' the definition.
%We'll see that they give the same notion.

The first idea, and the one leading to the definition given in \cite{HHS2}, is that not only the projections to the $\fontact U$ should be quasiconvex, but they should also determine the subspace. One ensures that this is the case by requiring that all realisation points of coordinates $(b_U\in\pi_U(Q))$ lie close to $Q$, where $Q$ is the hierarchically quasiconvex subspace. That also ensures that $Q$ is an HHS itself, see \cite[Proposition 5.6]{HHS2}.

A picture that one should keep in mind regarding this second property is that it is not satisfied by an "L" in $\mathbb{R}^2$, even though its projections to the two factors are convex. Rather, what one wants is a "full" square. This is related to the difference, in the cubical world, between $\ell_1$--isometric embeddings (the "L") as opposed to convex embeddings (the square).

% There is another idea, which is probably more intuitive. (Unfortunately, so far it seems less useful in applications and had note been recorded in the literature before; we do so in Lemma \ref{lem:hqconv_hierarchy} below.) Recall that given an HHS, there is a preferred family of quasigeodesics connecting pairs of points, which are called hierarchy paths, and that their defining property is that they project to unparametrized quasigeodesics in all $\fontact U$. It is then natural to just replace geodesics in the definition of the usual quasiconvexity with hierarchy paths, namely we want to say that $Q$ is hierarchically quasiconvex if all hierarchy paths joining points of $Q$ stay close to $Q$. It should be easy to see that this implies that the projections of $Q$ to the $\fontact U$ are quasiconvex, and we will argue that the additional property is satisfied too below, after we give formal definitions.

The second idea to complete the definition is probably more intuitive. Unfortunately, as of yet it has not been proven that this gives an equivalent notion. Recall that given an HHS, there is a preferred family of quasigeodesics connecting pairs of points, which are called hierarchy paths, and that their defining property is that they project to unparametrized quasigeodesics in all $\fontact U$. It is then natural to just replace geodesics in the definition of the usual quasiconvexity in hyperbolic spaces with hierarchy paths. Namely we want to say that $Q$ is hierarchically quasiconvex if all hierarchy paths joining points of $Q$ stay close to $Q$. It is easy to see that this implies that the projections of $Q$ to the $\fontact U$ are quasiconvex, but the additional property is not clear, and it would be interesting to known whether it is satisfied. It is however true (and easy to see) that hierarchy paths joining points on a hierarchically quasiconvex set stay close to it. Moreover, the hull of two points, as defined in the next subsection, coarsely coincides with the union of all hierarchy paths (with fixed, large enough constant) connecting them.

\begin{defn}[Hierarchical quasiconvexity]\label{defn:hierarchical_quasiconvexity}\cite[Definition 5.1]{HHS2}
Let $(\cuco X,\mathfrak S)$ be a hierarchically hyperbolic space.  
Then $\cuco Y\subseteq\cuco X$ is \emph{$k$--hierarchically quasiconvex}, for some $k\co[0,\infty)\to[0,\infty)$, if the following hold:
\begin{enumerate}
 \item For all $U\in\mathfrak S$, the projection $\pi_U(\cuco Y)$ is 
 a $k(0)$--quasiconvex subspace of the $\delta$--hyperbolic space $\fontact U$.
 \item For all $\kappa\geq0$ and $\kappa$-consistent tuples $\tup
 b\in\prod_{U\in\mathfrak S}2^{\fontact U}$ with $b_U\subseteq\pi_U(\cuco Y)$ for all $U\in\mathfrak S$, each point
 $x\in\cuco X$ for which $\dist_{
 U}(\pi_U(x),b_U)\leq\theta_e(\kappa)$ (where $\theta_e(\kappa)$ is as in Theorem~\ref{thm:realization}) satisfies $\dist(x,\cuco Y)\leq
 k(\kappa)$.
\end{enumerate}
\end{defn}

\begin{rem} Note that condition (2) in the above definition is 
    equivalent to:  For every $\kappa>0$ and every point 
    $x\in\cuco X$ satisfying 
    $\dist_{U}(\pi_U(x),\pi_{U}(\cuco Y))\leq\kappa$ 
    for all $U\in\mathfrak S$, has the property that 
    $\dist(x,\cuco Y)\leq k(\kappa)$.
\end{rem}

A very important property of hierarchically quasiconvex subspaces is that they admit a natural coarse retraction, which generalises the retraction onto standard product regions. The definition (see \cite[Definition 5.4]{HHS2}) is that a coarsely Lipschitz map 
$\gate_{\cuco Y}\colon\cuco X\to\cuco Y$ is called a \emph{gate map} 
if for each $x\in\cuco X$ it satisfies:  $\gate_{\cuco
Y}(x)$ is a point $y\in\cuco Y$ 
such that for all $V\in\mathfrak S$,
the set $\pi_V(y)$ (uniformly) coarsely coincides with the projection
of $\pi_V(x)$ to the $k(0)$--quasiconvex set $\pi_V(\cuco Y)$.

The uniqueness axiom implies that when such a map exists it 
is coarsely well-defined, and the existence is proven in \cite[Lemma 5.5]{HHS2}.

\subsection{Hulls and their cubulation}

Other important examples of hierarchically quasiconvex subspaces are hulls. In a $\delta$-hyperbolic space $X$, if one considers any subset $A\subseteq X$, then one can construct a ``quasiconvex hull'' simply by taking the union of all geodesics connecting pairs of points in the space. It is easy to show that this is a $2\delta$-quasiconvex subspace, whatever $A$ is. In an HHS, we can proceed similarly.

\begin{defn}[Hull of a set; \cite{HHS2}]\label{defn:hull}
 For each $A\subset \cuco X$ and $\theta\geq 0$, let the \emph{hull}, 
 $H_{\theta}(A)$,  be the set of all $p\in\cuco X$ so that, for each
 $W\in\mathfrak S$, the set $\pi_W(p)$ lies at distance at most
 $\theta$ from $\hull_{\fontact W}(A)$, the convex hull of $A$ in the 
 hyperbolic space $\fontact W$ (that is to say, the union of all geodesics in $\fontact W$ joining points of $A$). 
 Note that $A\subset H_{\theta}(A)$.
\end{defn}

It is proven in \cite[Lemma 6.2]{HHS2} that for each sufficiently large $\theta$ there
exists $\kappa\co\reals_+\to\reals_+$ so that for each $A\subset 
\cuco X$ the set $H_\theta(A)$ is
$\kappa$--hierarchically quasiconvex.

To connect a few things we have seen so far, the gate map onto the hull of two points $x,y$ coincides with taking the median $\median(x,y,\cdot)$.

Hulls of finite sets carry more structure than just being hierarchically quasiconvex. In a hyperbolic space, hulls of finitely many points are quasi-isometric to trees, with constants depending only on how many points one is considering. Trees are 1-dimensional CAT(0) cube complexes, and what happens in more general HHSs is that hulls of finitely many points are quasi-isometric to possibly higher dimensional CAT(0) cube complexes:

 \begin{thm}\label{thm:cubulated_hull}\cite[Theorem C]{HHS_quasiflats}
  Let $(\cuco X,\mathfrak S)$ be a hierarchically hyperbolic space and let $k\in\naturals$.  Then there exists $M_0$ so that for all $M\geq M_0$ there is a constant $C_1$ so that for any $A\subset\cuco X$ of cardinality $\leq k$, there is a $C_1$--quasimedian\footnote{CAT(0) cube complexes are endowed have a natural median map. Here by $C_1$--quasimedian we mean that the coarse median of 3 points is mapped within distance $C_1$ of the median of their images.}
  $(C_1,C_1)$--quasi-isometry $\mathfrak p_A\co\cuco Y\to H_\theta(A)$.
  
  Moreover, let $\mathcal U$ be the set of $U\in\mathfrak S$ so that $\dist_U(x,y)\geq M$ for some $x,y\in A$.  Then $\dimension\cuco Y$ is equal to the maximum cardinality of a set of pairwise-orthogonal elements of $\mathcal U$.
  
  Finally, there exist $0$--cubes $y_1,\ldots,y_{k'}\in\cuco Y$ so that $k'\leq k$ and $\cuco Y$ is equal to the convex hull in $\cuco Y$ of $\{y_1,\ldots,y_{k'}\}$.
 \end{thm}
 
 The theorem above is crucial for the proof of the quasiflats theorem in \cite{HHS_quasiflats}, and it definitely feels like it should have many more applications.

 \subsection{Factored spaces}
 
 In this section we discuss a construction where one cones-off subspaces of an HHS and obtains a new HHS. Examples of HHSs that are obtained using this construction from some other HHS are the ``main'' hyperbolic space $\fontact S$, which is obtained coning-off standard product regions in the corresponding $\cuco X$, and the Weil-Petersson metric on Teichm\"uller space, which is (coarsely) obtained coning-off Dehn twist flats in the corresponding mapping class group. The construction was devised in \cite{HHS_asdim}, with the purpose of having spaces that ``interpolate'' between a given HHS $\cuco X$ and the hyperbolic space $\fontact S$. In particular, this gives one way of doing induction arguments.
 
 Given a hierarchical space $\HHS X S$, we say $\mathfrak U\subseteq\mathfrak S$ is \emph{closed under nesting} if for all $U\in\mathfrak U$, if $V\in\mathfrak S-\mathfrak U$, then $V\not\nest U$.

\begin{defn}[Factored space]\label{defn:factored_space}\cite[Definition 2.1]{HHS_asdim}
Let $(\cuco X,\mathfrak S)$ be a hierarchically hyperbolic space. 
A \emph{factored space} $\coneoff{\cuco X}_{\mathfrak U}$ is constructed by
defining a new metric $\hat\dist$ on $\cuco X$ depending on a 
given subset $\mathfrak U \subset\mathfrak S$ which is closed under nesting. 
First, for each $U\in\mathfrak U$, for each pair $x,y\in\cuco X$ for 
which there exists $e\in {\bf E}_{U}$ such that $x,y\in {\bf 
F}_{U}\times \{e\}$, we set 
$\dist'(x,y)=\min\{1,\dist(x,y)\}$.  
For any pair $x,y\in\cuco X$ for which there does not exists such an 
$e$ we set $\dist'(x,y)=\dist(x,y)$.  We now
define the distance $\hat\dist$ on $\coneoff{\cuco X}_{\mathfrak U}$.  Given a
sequence $x_0,x_1,\ldots,x_k\in\coneoff{\cuco X}_{\mathfrak U}$, define its length
to be $\sum_{i=1}^{k-1}\dist'(x_i,x_{i+1})$.  Given
$x,x'\in\coneoff{\cuco X}_{\mathfrak U}$, let $\hat\dist(x,x')$ be the infimum of
the lengths of such sequences $x=x_0,\ldots,x_k=x'$.  
\end{defn}

 It is proven in \cite[Proposition 2.4]{HHS_asdim} that factored spaces are HHSs themselves, with the ``obvious'' substructure of the original HHS structure.
 
 I believe that factored spaces will play an important role in studying quasi-isometries between hierarchically hyperbolic spaces in view of \cite[Corollary 6.3]{HHS_quasiflats}, which says that a quasi-isometry between HHSs induces a quasi-isometry between certain factored spaces. One can then take full advantage of the HHS machinery to iterate this procedure and get more information about the quasi-isometry under examination. This strategy should work in various sub-classes of right-angled Artin and Coxeter groups (one example is given in \cite{HHS_quasiflats}).

\subsection{Boundary}

HHSs admit a boundary, defined in \cite{HHS_boundary}. Here is the idea behind the definition. We have to understand how one goes to infinity in an HHS, and more precisely we want to understand how one goes to infinity along a hierarchy ray, because we like hierarchy paths. The hierarchy ray will project to any given $\fontact U$ close to either a geodesic segment or a geodesic ray. In fact, a consequence of large links is (the hopefully intuitive fact) that there must be some $\fontact U$ so that the hierarchy ray has unbounded projection there (see the useful ``passing up'' lemma \cite[Lemma~2.5]{HHS2}). Now, one can consider the non-empty set of all $\fontact U$ where the hierarchy ray has unbounded projection, and it should come to no surprise that they must be pairwise orthogonal, because in all other cases there are constraints coming from the consistency inequalities on the projections of points of $\cuco X$. As a last thing, you might want to measure how fast you are moving asymptotically in each of the $\fontact U$ and record the ratios of the various speeds. These considerations lead to the following definition. 

\begin{defn}\label{defn:boundary_point}\cite[Definitions 2.2-2.3]{HHS_boundary}
A \emph{support set} $\overline S\subset\mathfrak S$ is a set with $S_i\orth S_j$ for all $S_i,S_j\in\overline S$.  Given a support set $\overline S$, a \emph{boundary point} with \emph{support} $\overline S$ is a formal sum $p=\sum_{S\in\overline S}a^p_Sp_S$, where each $p_S\in\boundary\fontact S$, and $a^p_S>0$, and $\sum_{S\in\overline S}a^p_S=1$.  (Such sums are necessarily finite.)
%We denote the support $\overline S$ of $p$ by $\support(p)$.

The \emph{boundary} $\boundary(\cuco X,\mathfrak S)$ of $(\cuco X,\mathfrak S)$ is the set of boundary points.
\end{defn}

The boundary also comes with a topology, which is unfortunately very complicated to define. But fortunately, for at least some applications one can just use a few of its properties and never work with the actual definition (this is the case for the rank rigidity theorems \cite[Theorems 9.13,914]{HHS_boundary}). The main property is that if an HHS is proper as a metric space then its boundary is compact by \cite[Theorem 3.4]{HHS_boundary}. Another good property is that hierarchically quasiconvex subspaces (and not only those) have well-defined boundary extensions.

% For applications and more discussion, see the introduction of \cite{HHS_boundary, Mousley1,Mousley2}.

\subsection{Modifying the HHS structure}

There are various way to modify an HHS structure, and this is useful to perform certain constructions. For example, there is a combination theorem for trees of HHSs which, unsurprisingly, requires the HHS structures of the various edge and vertex groups to be ``compatible''. In concrete cases, one might be starting with a tree of HHSs that does not satisfy the compatibility conditions on the nose, but does after suitably adjusting the various HHS structures. This is the case even for very simple examples like the amalgamated product of two mapping class groups over maximal virtually cyclic subgroups containing a pseudo-Anosov.

A typical way of modifying the HHS structure is to cone-off a collection of quasi-convex subspaces of one of the $\fontact U$, and then add those quasi-convex subspaces as new hyperbolic spaces. This usually amounts to considering certain (hierarchically quasiconvex and hyperbolic) subspaces of $\cuco X$ as product regions with a trivial factor. This construction is used to set things up for the small cancellation constructions in \cite{HHS_asdim}, and is explored in much more depth in \cite{Spriano}. In the latter paper the author studies very general families of quasiconvex subspaces of $\fontact S$ (called factor systems) so that one can perform a version of the cone-off-and-add-separately construction mentioned above. This is needed to set things up before using the combination theorem in natural examples, like the one mentioned above.

In another direction, one may wonder whether there exists a minimal HHS structure, and in particular one might want to reverse the procedure above. This is explored in \cite{HHS_stable}.

\bibliographystyle{alpha}
\bibliography{document}
\end{document}